\documentclass[12 pt]{article}
\usepackage{amsmath}
\usepackage{amssymb}
\usepackage{amsthm}
\usepackage{graphicx}
\usepackage{subfig}
\usepackage{float}
\usepackage{setspace}

\newtheorem{theorem}{Theorem}
\newtheorem{lemma}{Lemma}

\newtheorem{corollary}{Corollary}
\newtheorem{definition}{Definition}
\newtheorem{remark}{Remark}

\title{Relative Equilibria of the $(1+N)$-Vortex Problem}
\author{Anna M. Barry$^1$\\Department of Mathematics and Statistics\\Center for BioDynamics\\Boston University \and Glen R. Hall\\Department of Mathematics and Statistics\\Boston University \and C. Eugene Wayne\\Department of Mathematics and Statistics\\Center for BioDynamics\\Boston University}

\begin{document}

\maketitle

\footnote[1]{Corresponding author: annab@math.bu.edu}

\newpage

\singlespacing

\begin{abstract}
We examine existence and stability of relative equilibria of the $n$-vortex problem specialized to the case where $N$ vortices have small and equal circulation and one vortex has large circulation.  As the small circulation tends to zero, the weak vortices tend to a circle centered on the strong vortex.  A special potential function of this limiting problem can be used to characterize orbits and stability.  Whenever a critical point of this function is nondegenerate, we prove that the orbit can be continued via the Implicit Function Theorem, and its linear stability is determined by the eigenvalues of the Hessian matrix of the potential.  For $N\geq 3$ there are at least three distinct families of critical points associated to the limiting problem.  Assuming nondegeneracy, one of these families continues to a linearly stable class of relative equilibria with small and large circulation of the same sign.  This class becomes unstable as the small circulation passes through zero and changes sign.  Another family of critical points which is always nondegenerate continues to a configuration with small vortices arranged in an $N$-gon about the strong central vortex.  This class of relative equilibria is linearly unstable regardless of the sign of the small circulation when $N\geq 4$.  Numerical results suggest that the third family of critical points of the limiting problem also continues to a linearly unstable class of solutions of the full problem independent of the sign of the small circulation.  Thus there is evidence that linearly stable relative equilibria exist when the large and small circulation strengths are of the same sign, but that no such solutions exist when they have opposite signs.  The results of this paper are in contrast to those of the analogous celestial mechanics problem, for which the $N$-gon is the only relative equilibrium for $N$ sufficiently large, and is linearly stable if and only if $N\geq 7$.
\end{abstract}

\newpage

\section{Introduction}\label{sect:1}

\indent A well known problem in fluid mechanics is that of $n$ point vortices moving in the plane, which may be regarded as an approximation to more general vortex motion.  This is commonly known as the \textit{$n$-vortex problem}.  Of particular interest are solutions that appear fixed when viewed in a uniformly rotating or translating frame.  These are called \textit{relative equilibria} and have been studied by many authors, see \cite{MR1831715}, \cite{MR1740937}, \cite{MR1332658}, \cite{MR1411341}, \cite{MR648066}, \cite{MR1753026} and \cite{MR2337012}.\newline
\indent The focus of this paper is on relative equilibria of the planar $n$-vortex problem with $n=N+1$ where one vortex has large circulation and $N$ have small, equal circulation.  We show that the weak vortices limit to a single circle about the strong central vortex as the small circulation tends to zero, and that the resulting configuration must be a critical point of a particular potential function depending only on the relative positions of the ``small'' vortices on this limiting circle.  Conversely, we find that nondegenerate critical points of this potential continue to relative equilibria with nonzero circulation, and the linear stability of these solutions is determined by the type of critical point of the limit potential.  \newline
\indent One case that has been carefully studied has the small vortices in an $N$-gon around the large vortex in the center.  This family is a relative equilibrium of the problem regardless of the relative strengths of the central and surrounding vortices.  Its stability has been studied in detail by Cabral and Schmidt in \cite{MR1740937}.  In particular, they showed that the $N$-gon around a central vortex is unstable when the circulation of the surrounding vortices is sufficiently small.  A new proof of this is presented here and these techniques imply the existence of two more families of critical points of the limit potential.  One of these other families continues to a linearly stable class of relative equilibria of the full problem whenever it is nondegenerate.  Members of this class of configurations have weak and strong vortices of the same sign.\newline
\indent The methods used here are parallel to those used to study relative equilibria of the $(1+N)$-body problem with Newtonian inverse square force of attraction in celestial mechanics.  The results for vortices are, surprisingly, very different.  For the celestial mechanics problem, the relative equilibrium formed by an $N$-gon of small masses about a large central mass is stable for $N\geq 7$ and is the only relative equilibrium for $N$ sufficiently large, see \cite{hall}, \cite{MR1472896}, \cite{MR1262722}, \cite{MR1816912}, \cite{MR1310189} and \cite{MR2104897}.  Another difference is that the circulation parameter of the vortex problem, which is analogous to mass, can be positive or negative.  We will see that there can be a dramatic bifurcation as this parameter passes through zero.\newline
\indent In the next section of the paper we carefully define the problem and build the important potential function.  In Section \ref{sect:3} we carry out the proof that nondegenerate critical points of this potential can be continued to relative equilibria with $N$ vortices having small, nonzero circulation and one having large circulation.  Moreover, we show that the linear stability of such a continuation is determined by properties of the potential.  Section \ref{sect:4} consists of applications of Section \ref{sect:3} results to the cases $N=2$, $3$ and $4$.  In Section \ref{sect:5} we discuss further results for large $N$ and do some numerical investigation.  In particular, we prove that when $N\geq 4$ the $N$-gon consisting of $N$ small vortices surrounding a strong central vortex is linearly unstable.  We use this to show that when all critical points are nondegenerate there must exist at least two other families of relative equilibria of the full problem.  We argue that one of these families is made up of linearly stable orbits when the central and surrounding vortices are of the same sign.  We use numerics to study the candidates for the three families and carry out stability calculations.

\section{Relative Equilibria}\label{sect:2}

Consider an inviscid, incompressible fluid in the plane with vorticity distribution described by a finite sum of Dirac delta functions.  Note that this ``function'' is zero except at the locations of the point vortices.  Associated to each of these point vortices is an important quantity $\Gamma(t)$, called the \textit{circulation} of the vortex.  Assuming the flow is acted on only by conservative forces, we may conclude that $\Gamma(t)\equiv \Gamma$, a constant.  We therefore treat $\Gamma$ as a parameter of the problem.\newline
\indent With this in mind, let $q_i(t)=(x_i(t),y_i(t))\in\mathbb{R}^2$, $i=1,...,n$ be the positions of $n$ point vortices in the plane with associated circulations $\Gamma_i$.  Then the equations of motion are described by
\begin{align}
\Gamma_i\dot{x_i}&=\frac{\partial H}{\partial y_i}\\
\Gamma_i\dot{y_i}&=-\frac{\partial H}{\partial x_i}
\end{align}
where the Hamiltonian function $H$ is
\begin{equation}
H(q_1,...,q_n)=-\sum_{i,j=1,i<j}^n \Gamma_i\Gamma_j\log|q_i-q_j|.
\end{equation}
For details, see \cite{MR1831715}.\newline
\indent For $n=2,3$ the system is completely integrable, and all possible relative equilibria have been discovered.  Proof and examples can be found in \cite{MR1831715}.  Later, we will discuss how these orbits appear in the specialized problem of one strong and $N$ weak vortices. For $n>3$ the problem is nonintegrable, and hence more difficult.  For proof of the nonintegrability of the four-vortex problem, see \cite{MR564329} or \cite{MR987772}.\newline
\indent We now specialize the general problem to our desired setting.  Consider $N+1$ vortices $q_0,q_1,...,q_N\in\mathbb{R}^2$ with strengths $\Gamma_0=1$, $\Gamma_i=\varepsilon$, $i=1,...,N$.  The equations can be written in the form
\begin{equation}\label{eqn:1}
\dot{q_j}=\sum_{i=0,i\neq j}^N \Gamma_i\frac{(q_j-q_i)^\bot}{|q_j-q_i|^2},\;\; j=0,...,N
\end{equation}
where $(x,y)^\bot=(-y,x)$.  Next we define what it means to be a relative equilibrium of this specific problem and explore some properties of solutions.

\begin{definition}\label{defn:1} Let $\{\varepsilon_k\}_{k=1}^\infty$ be a sequence of real numbers such that $\varepsilon_k\rightarrow 0$ as $k\rightarrow\infty$ and let $q_0^k,...,q_N^k$ be a sequence of configurations which are relative equilibria of the $(N+1)$-vortex problem for each $k$ with corresponding circulations $\Gamma_0^k=1,\Gamma_j^k=\varepsilon_k$, $j=1,...,N$.  A \textup{relative equilibrium of the $(1+N)$-vortex problem} is a configuration $q_0,...,q_N$ such that $q_j^k\rightarrow q_j$ as $k\rightarrow\infty$ for each $j=0,...,N$.\end{definition}

\begin{remark}\label{rmk:1} Both positive and negative circulation $\varepsilon$ make physical sense.  The existence and continuation results that follow hold for $\varepsilon$ of either sign.  However, linear stability depends on the sign of $\varepsilon$, which we point out below.\end{remark}

\begin{remark} \label{rmk:2} To avoid confusion, we remark that the phrase \textup{$(N+1)$-vortex problem} is meant to refer to the original equations (\ref{eqn:1}) of the $n$-vortex problem with $n=N+1$.  On the other hand, the phrase \textup{$(1+N)$-vortex problem} refers to the limiting problem in which there is one strong vortex and $N$ vortices with zero limiting circulation.  This is consistent with the celestial mechanics literature.\end{remark}

Definition \ref{defn:1} does not rule out the possibility that two or more vortices collide in the limit.  However, we wish to exclude this situation from our calculations and so we require that all vortices are bounded away from each other by some constant $m>0$.  We also ignore configurations which become unbounded as $k\rightarrow\infty$, and so we assume that all vortices lie inside a circle of radius $M$ about the origin.  Both $M$ and $m$ are assumed to be independent of $k$ and $\varepsilon_k$, but may depend on $N$.  Under these assumptions we can prove a lemma which asserts that relative equilibria of the $(1+N)$-vortex problem have small vortices on a circle of fixed radius about the large vortex, which is also fixed.  Since the center of vorticity is an invariant of the problem we take it to be the origin for simplicity, i.e. we choose $q_0+\varepsilon(q_1+...+q_N)=0$.  This will ensure that the large vortex is near the origin.  Moreover, by choosing an appropriate rotation rate the small vortices will limit to the unit circle.\newline
\indent The next two results have counterparts in celestial mechanics, see \cite{hall} or \cite{MR1262722}.

\begin{lemma} \label{lma:1} In a rigidly rotating frame with rotation rate $\omega=1$, all relative equilibria, $q_0^k,...,q_N^k$, which converge to a relative equilibrium $q_0,...,q_N$ of the $(1+N)$-vortex problem satisfy $|q_0^k|=\mathcal{O}(\varepsilon_k)$, $|q_j^k|^2-1=\mathcal{O}(\varepsilon_k)$, $j=1,...,N$.
\end{lemma}

\begin{proof} We will drop the superscript $k$ whenever it is unambiguous to do so, i.e. when $\varepsilon\neq 0$ is clear. Let $q_j^k=(x_j^k,y_j^k)\in\mathbb{R}^2$, $z_j^k=x_j^k+iy_j^k$, and $\xi_j^k=e^{i\omega t}z_j^k$.  Suppose $\xi^k=(\xi_0^k,...,\xi_N^k)$ is a fixed point of the system (\ref{eqn:1}) written in rotating coordinates.  That is,
\begin{align}
0&=-i\omega\xi_0+\varepsilon\sum_{i=1}^N \frac{1}{\xi_0^*-\xi_i^*}\\
0&=-i\omega \xi_j+\frac{1}{\xi_j^*-\xi_0^*}+\varepsilon\sum_{i=1,i\neq j}^N \frac{1}{\xi_j^*-\xi_i^*},\;\; j=1,...,N.
\end{align}
For fixed $\omega=1$,
\begin{align}
i\xi_0&=\varepsilon\sum_{i=1}^N \frac{1}{\xi_0^*-\xi_i^*}\\
i\xi_j&=\frac{1}{\xi_j^*+\varepsilon(\xi_1^*+...+\xi_N^*)}+\varepsilon\sum_{i=1,i\neq j}^N \frac{1}{\xi_j^*-\xi_i^*}
\end{align}
and so we see that
\begin{equation}
|\xi_0|\leq \varepsilon\frac{N}{m}
\end{equation}
which proves the first assertion of the lemma.  From this and the assumption that the distances between vortices is bounded below we have
\begin{equation}
|\xi_j|=|\xi_0-\xi_j-\xi_0|\geq |\xi_0-\xi_j|-|\xi_0|\geq m-\varepsilon \frac{N}{m}
\end{equation}
Therefore,
\begin{align}
|\xi_j|^2&\leq |1-\frac{\varepsilon}{\xi_j^*}\left(\xi_1^*+...+\xi_N^*\right)+\mathcal{O}(\varepsilon^2)|+\frac{\varepsilon M(N-1)}{m}\\
&\leq 1+\varepsilon \left(\frac{M N}{m-\varepsilon\frac{N}{m}}+\frac{M(N-1)}{m}\right)+\mathcal{O}(\varepsilon^2)\\
&= 1+\varepsilon \left(\frac{M N}{m}+\frac{M(N-1)}{m}\right)+\mathcal{O}(\varepsilon^2).
\end{align}
Similarly, we find
\begin{equation}
|\xi_j|^2\geq 1-\frac{\varepsilon N M}{m}-\frac{\varepsilon M(N-1)}{m}+\mathcal{O}(\varepsilon^2)
\end{equation}
which concludes the proof.\end{proof}

Since any periodic relative equilibrium must rotate rigidly about its center of vorticity (assumed to be zero), a necessary condition for a formation to be in relative equilibrium is $q_j\cdot \dot{q_j}=0$, $j=1,...,n$.  Using this idea, we prove the following lemma.

\begin{lemma}\label{lma:2} Let $q=(q_0,...,q_N)$ be a relative equilibrium of the $(1+N)$-vortex problem and let $(r,\theta)=(r_1,...,r_N,\theta_1,...\theta_N)$ be the representation of $(q_1,...,q_N)$ in polar coordinates.  Then $\theta$ is a critical point of the potential function
\begin{equation}\label{eqn:2}
V(\theta)=-\sum_{i<j}\left(\cos(\theta_i-\theta_j)+\frac{1}{2}\log(2-2\cos(\theta_i-\theta_j))\right).
\end{equation}
\end{lemma}

\begin{proof}  Let $q^\varepsilon$ be a sequence of relative equilibria of the full $(N+1)$-vortex problem which converges to $q$ as $\varepsilon\rightarrow 0$. Since $q_j^\varepsilon\cdot \dot{q}_j^\varepsilon=0$, after suppressing the $\varepsilon$ dependence we have, for $j=1,...,N$,
\begin{align}
0&=q_j\cdot\left(\frac{q_j^\bot-q_0^\bot}{|q_j-q_0|^2}+\sum_{i=1,i\neq j}^N\frac{\varepsilon(q_j^\bot-q_i^\bot)}{|q_j-q_i|^2}\right)\\
&=q_j\cdot \left(\frac{q_j^\bot+\varepsilon(q_1+...+q_N)^\bot}{|q_j+\varepsilon(q_1+...+q_N)|^2}+\sum_{i=1,i\neq j}^N\frac{\varepsilon(q_j^\bot-q_i^\bot)}{|q_j-q_i|^2}\right)\\
&=\left(\frac{\varepsilon q_j\cdot(q_1+...+q_N)^\bot}{|q_j+\varepsilon(q_1+...+q_N)|^2}-\sum_{i=1,i\neq j}^N\frac{\varepsilon(q_j\cdot q_i^\bot)}{|q_j-q_i|^2}\right)\\
&=\varepsilon\sum_{i=1,i\neq j}^N q_j\cdot q_i^\bot\left(\frac{1}{|q_j+\varepsilon(q_1+...+q_N)|^2}-\frac{1}{|q_j-q_i|^2}\right).
\end{align}
In polar coordinates, an expansion of the above expression in powers of $\varepsilon$ is given by
\begin{equation}\label{eqn:3}
0=\varepsilon\sum_{i\neq j}^N r_j r_i \sin(\theta_j-\theta_i)\left(\frac{1}{r_j^2}-\frac{1}{r_j^2+r_i^2-2r_j r_i\cos(\theta_j-\theta_i)}\right)+\mathcal{O}(\varepsilon^2).
\end{equation}
Due to Lemma \ref{lma:1}, dividing both sides by $\varepsilon$ and taking the limit as $\varepsilon\rightarrow 0$ yields
\begin{equation}
0=\sum_{i\neq j}^N \left(\sin(\theta_j-\theta_i)-\frac{\sin(\theta_j-\theta_i)}{2-2\cos(\theta_j-\theta_i)}\right).
\end{equation}
Therefore, $\theta$ is a critical point of
\begin{equation}
V(\theta)=-\sum_{i<j}\left(\cos(\theta_i-\theta_j)+\frac{1}{2}\log(2-2\cos(\theta_i-\theta_j))\right).
\end{equation}
\end{proof}

\begin{remark}\label{rmk:3} In what follows, the Hessian matrix $V_{\theta\theta}$ will play an important role.  Due to the rotational symmetry of the problem, this matrix will always have at least one zero eigenvalue with corresponding eigenvector $v_0:=(1,1,...,1)^T$.  Therefore, as pointed out in \cite{MR1472896}, it is appropriate to define a critical point, $\theta$, of $V$ to be \textit{nondegenerate} provided $V_{\theta\theta}$ has only one zero eigenvalue.  Further, we say that $\theta$ is a \textit{nondegenerate local minimum} (respectively, maximum) of the potential $V$ if $V_{\theta\theta}$ is positive (negative) semidefinite with a one-dimensional null space.\end{remark}

\section{Linear Stability}\label{sect:3}

\indent Due to the many symmetries and integrals present in our system, the usual stability theory for ordinary differential equations does not apply.  Here, we modify the definitions of nondegeneracy and linear stability to account for these degeneracies.  \newline
\indent First of all, it is impossible for a fixed point (or periodic orbit) of the system (\ref{eqn:1}) to be hyperbolic because there are always at least four zero eigenvalues associated to the linearization.  Two of these eigenvalues come from the invariance of the center of vorticity, which we have assumed is at the origin. By defining $q_0:=-\varepsilon(q_1+...+q_N)$, the equations for the components of $q_0$ are satisfied once $q_1,...,q_N$ have been specified.  Therefore, we ignore these two equations and focus only on the remaining $2N$ equations, hence removing two zero eigenvalues from the system.  Written in polar coordinates, the reduced system becomes, for $j=1,...,N$,
\begin{align}
\frac{d}{dt}(r_j^2)&=
\varepsilon\sum_{i\neq j}^N r_j r_i \sin(\theta_j-\theta_i)\left(\frac{1}{r_j^2}-\frac{1}{r_j^2+r_i^2-2r_j r_i\cos(\theta_j-\theta_i)}\right)+\mathcal{O}(\varepsilon^2) \notag\\
&=:\varepsilon F_j(r,\theta,\varepsilon)\label{eqn:4}\\
\frac{d\theta_j}{dt}&=\frac{1}{r_j^2}+\varepsilon\sum_{i\neq j} \frac{r_i^2(r_j\cos(2(\theta_i-\theta_j))-r_i\cos(\theta_i-\theta_j))}{r_j^3(r_j^2+r_i^2-2 r_i r_j \cos(\theta_i-\theta_j))}+\mathcal{O}(\varepsilon^2).\notag\\
&=:G_j(r,\theta,\varepsilon).\label{eqn:5}
\end{align}
The remaining two zero eigenvalues can be readily identified once we linearize this system about a relative equilibrium of the problem.  Consider a fixed point of (\ref{eqn:4})-(\ref{eqn:5}) which is a member of a sequence of relative equilibria converging to a relative equilibrium of the $(1+N)$-vortex problem.  Then upon applying Lemma \ref{lma:1} for $\varepsilon$ sufficiently small, the matrix $M$ of the linearized problem about this fixed point is made up of four $N\times N$ blocks:

\begin{equation}\label{eqn:6}
M=
\left( {\begin{array}{cc}
 -\varepsilon A +\mathcal{O}(\varepsilon^2) & \varepsilon V_{\theta\theta}(\phi)+\mathcal{O}(\varepsilon^2) \\
 -2 I+\mathcal{O}(\varepsilon) & \varepsilon A+\mathcal{O}(\varepsilon^2)  \\
 \end{array} } \right)
\end{equation}

\noindent where $A=(a_{ij})$ satisfies
\begin{align}
a_{ij}&=\sin(\phi_j-\phi_i)\;\; i\neq j\\
a_{ii}&=\sum_{j\neq i} \sin(\phi_i-\phi_j).
\end{align}

From the form of $A$ and Remark \ref{rmk:3}, it is easy to see that $M (0,v_0)^T=\mathcal{O}(\varepsilon^2)$.  This computation suggests that one zero eigenvalue may be related to the rotational symmetry of the problem.  This is true, and we leave the details as an exercise.  We show in the proof of Theorem \ref{thm:1} that we may ignore this eigenvalue by focusing on the complementary eigenspace to the span of $v_0$.\newline
\indent Due to the Hamiltonian nature of the problem, this eigenvalue must be paired with a second zero eigenvalue. Physically, this second eigenvalue is born from what can be thought of as the scaling symmetry of the problem: associated to any relative equilibrium is a family of relative equilibria which are related to the first by a scaling and have different rotation frequencies.  One may deal with this degeneracy by specifying a rotation frequency, as we did in the proof of Lemma \ref{lma:1}.  \newline
\indent With this in mind, it is reasonable to define a relative equilibrium of the $(N+1)$-vortex problem to be \textit{nondegenerate} if the associated linearization has exactly four zero eigenvalues (or, once the equations for $q_0$ are removed as described, exactly two zero eigenvalues).\newline
\indent Aside from the zero eigenvalues,  we remind the reader that for Hamiltonian systems asymptotic stability is impossible, since all eigenvalues come in pairs with opposite signs.  Therefore, one says instead that an equilibrium is linearly stable if all eigenvalues are purely imaginary.  In the proof of the main theorem of this section, we leave the two zero eigenvalues associated to scaling and rotational symmetry ``as is'' and call an equilibrium \textit{linearly stable} if all remaining eigenvalues are purely imaginary and nonzero.\newline
\indent The first theorem of this section asserts the existence of a convergent $\varepsilon$-dependent family of relative equilibria whenever there is a nondegenerate critical point of the potential.

\begin{theorem} \label{thm:1} Let $\phi=(\phi_1,...,\phi_N)$ be a nondegenerate critical point of the potential $V$.  Then for $\rho=(1,1,...,1)$, the configuration $(\rho,\phi)$ is a relative equilibrium of the $(1+N)$-vortex problem, that is, there exists a sequence of relative equilibria of the $(N+1)$-vortex problem which converges to $(\rho,\phi)$ as $\varepsilon\rightarrow 0$.\end{theorem}

\begin{proof} The proof relies on two applications of the Implicit Function Theorem. Define $F(r,\theta,\varepsilon):=(F_1(r,\theta,\varepsilon),...,F_N(r,\theta,\varepsilon))$ and $G(r,\theta,\varepsilon):=(G_1(r,\theta,\varepsilon),...,G_N(r,\theta,\varepsilon))$ where $F_j$ and $G_j$ are is in (\ref{eqn:4})-(\ref{eqn:5}).
For any $\theta\in\mathbb{R}^n$, set $\theta=\hat{\theta}+\theta_{\mathrm{null}}$ where $\hat{\theta}\in\mathrm{span}\{v_0\}^\bot$ and $\theta_{\mathrm{null}}\in\mathrm{span}\{v_0\}$, $v_0$ as in Remark \ref{rmk:3}.  Then $F(r,\theta,\varepsilon)=F(r,\hat{\theta},\varepsilon)$ and $D_\theta F(r,\theta,\varepsilon)=D_\theta F(r,\hat{\theta},\varepsilon)$ because the equations depend only on the differences between the angular components.  Therefore, for the purposes of applying the Implicit Function Theorem, we restrict $F$ and $DF$ to span$\{v_0\}^\bot$.\newline
\indent With this in mind, set $\phi=\hat{\phi}+\phi_{\mathrm{null}}$.  Since $F(\rho,\hat{\phi},0)=0$ and the restriction of $D_\theta F(\rho,\hat{\phi},0)=V_{\theta\theta}(\hat{\phi})$ to this subspace is invertible, we obtain open sets $\mathcal{U},\mathcal{V}$ and $\mathcal{W}$ containing $r=\rho,$ $\theta=\hat{\phi}$ and $\varepsilon=0$, respectively, and a continuously differentiable function $f:\mathcal{U}\times\mathcal{W}\rightarrow \mathcal{V}$ such that $F(r,f(r,\varepsilon),\varepsilon)\equiv 0$ whenever $(r,\varepsilon)\in\mathcal{U}\times\mathcal{W}$.  This implies that there exists a sequence of solutions of the $(N+1)$-vortex problem which converges to $(\rho,\hat{\phi})$ as $\varepsilon\rightarrow 0$.  However, at this point it is not clear that these solutions are relative equilibria of the problem, since it is possible that different vortices have different rotation rates.\newline
\indent To finish the proof, we apply the Implicit Function Theorem to $G$.  Set $\theta=f(r,\varepsilon)$ with $(r,\varepsilon)$ restricted to $\mathcal{U}\times\mathcal{W}$, as determined in the previous paragraph.  Then $G(\rho,f(\rho,0),0)=G(\rho,\hat{\phi},0)=1$ and
\begin{align}
D_r(G(r,f(r,\varepsilon),\varepsilon))\big|_{(r,\varepsilon)=(\rho,0)}&=D_rG(\rho,\hat{\phi},0)+\frac{\partial f}{\partial r} (\rho,0)D_\theta G(\rho,\hat{\phi},0)\\
&=-2 I.
\end{align}
Thus there exist open sets $\mathcal{U}_1\subset \mathcal{U}$, $\mathcal{W}_1\subset\mathcal{W}$ containing $r=\rho$ and $\varepsilon=0$, and a continuous function $g:\mathcal{W}_1\rightarrow\mathcal{U}_1$ such that $G(g(\varepsilon),f(g(\varepsilon),\varepsilon),\varepsilon)\equiv 1$ whenever $\varepsilon\in\mathcal{W}_1$.  It is also true that $F(g(\varepsilon),f(g(\varepsilon),\varepsilon),\varepsilon)=0$ for $\varepsilon\in\mathcal{W}_1$ and so the proof is complete.
\end{proof}

The symplectic structure associated to any Hamiltonian system provides a rich theory and many tools for analysis.  Before stating the main theorem of this section, we list one definition and result that will be used in the proof.  However, we will not go into the details of the theory here.  For further discussion and proof of the lemma, see \cite{MR1472896}.

\begin{definition}\label{defn:2} Let $J$ be the $2N\times 2N$ block matrix given by
\[
J=
\left( {\begin{array}{cc}
 0 & -I \\
 I & 0 \\
 \end{array} } \right).
\]
The \textup{skew inner product} of two vectors $v,w\in\mathbb{R}^{2N}$ is defined to be
\begin{equation}
\Omega(v,w)=v^T J w.
\end{equation}\end{definition}

\begin{lemma}\label{lma:3} Suppose each eigenvector $v$ of $M$ satisfies the inequality $\Omega(v,v^*)\neq 0$.  Then all eigenvalues of $M$ are purely imaginary.\end{lemma}

\begin{theorem} \label{thm:2} Suppose $(\rho^\varepsilon,\phi^\varepsilon)=(\rho_1^\varepsilon,...,\rho_N^\varepsilon,\phi_1^\varepsilon,...,\phi_N^\varepsilon)$ is a sequence of relative equilibria of the $(N+1)$-vortex problem which converges to a relative equilibrium $(\rho,\phi)=(1,...,1,\phi_1,...,\phi_N)$ of the $(1+N)$-vortex problem as $\varepsilon\rightarrow 0$, where $\phi$ is a nondegenerate critical point of the potential $V$ defined by Equation (\ref{eqn:3}). Then $(\rho^\varepsilon,\phi^\varepsilon)$ is nondegenerate for $\varepsilon\neq 0$ sufficiently small. Moreover, in this case the configuration $(\rho^\varepsilon,\phi^\varepsilon)$ for $\varepsilon>0$ is linearly stable if and only if $\phi$ is a local minimum of $V$.  Likewise, for $\varepsilon<0$ it is linearly stable if and only if $\phi$ is a local maximum of $V$.\end{theorem}

\begin{proof} Much of the following proof is modeled on analogous arguments made by Moeckel in \cite{MR1262722} for the corresponding celestial mechanics problem. Let $M$ denote the linearization matrix (\ref{eqn:5}) about the relative equilibrium $(\rho^\varepsilon,\phi^\varepsilon)$.  Throughout the proof we suppress the explicit $\varepsilon$-dependence when there is no ambiguity.  Let $\lambda\in \mathbb{C}$ be nonzero and consider the matrix $M-\lambda I_{2N\times 2N}$:

\[
M-\lambda I=
\left( {\begin{array}{cc}
 -\lambda I-\varepsilon A+\mathcal{O}(\varepsilon^2) & \varepsilon V_{\theta\theta}+\mathcal{O}(\varepsilon^2) \\
 -2 I+\mathcal{O}(\varepsilon) & -\lambda I+\varepsilon A+\mathcal{O}(\varepsilon^2)  \\
 \end{array} } \right).
\]
Since the matrix $V_{\theta\theta}(\phi)$ is assumed to be nondegenerate, it has exactly one zero eigenvalue.  This, coupled with the discussion prior to Theorem \ref{thm:1} implies that $M$ has at least two zero eigenvalues.  Choose $0<c_1<1$ small enough that $|\zeta|> c_1$ for any nonzero eigenvalue, $\zeta$, of $V_{\theta\theta}(\phi)$.  We now consider $\det(M-\lambda I)$ restricted to the region $c_1 \sqrt{\varepsilon}\leq |\lambda|\leq \varepsilon^{1/4}$.  An application of Lemma \ref{lma:4} yields
\begin{align}
\det(M-\lambda I)&=\det(-\lambda I+\mathcal{O}(\varepsilon^2))\det(-\lambda I+\mathcal{O}(\varepsilon)\notag\\
&+(-2 I+\mathcal{O}(\varepsilon))(-\lambda I-\mathcal{O}(\varepsilon))^{-1}(\varepsilon V_{\theta\theta}+\mathcal{O}(\varepsilon^2)))\notag\\
&=\det(-\lambda I+\mathcal{O}(\varepsilon))\det(-\lambda I-\frac{2\varepsilon}{\lambda} V_{\theta\theta}+\mathcal{O}(\varepsilon))\\
&=\det(I+\mathcal{O}(\sqrt{\varepsilon}))\det(\lambda^2 I+2\varepsilon V_{\theta\theta}+\mathcal{O}(\varepsilon^{5/4}))\notag\\
&=(1+\mathcal{O}(\sqrt{\varepsilon}))\det(\lambda^2I +2\varepsilon V_{\theta\theta}+\mathcal{O}(\varepsilon^{5/4})).\notag
\end{align}
The only way that the right side can be zero in this region is if $\lambda =\mathcal{O}(\sqrt{\varepsilon})$.  Set $\lambda(\varepsilon)=\sqrt{\varepsilon}\gamma(\varepsilon)$ and define $\lim_{\varepsilon\rightarrow 0}\gamma(\varepsilon)=:\gamma_0$.  Then
\begin{equation}
\det(\varepsilon \gamma^2 I+2\varepsilon V_{\theta\theta}+\mathcal{O}(\varepsilon^{3/2}))=\varepsilon^N \det(\gamma^2 I+2V_{\theta\theta}+\mathcal{O}(\sqrt{\varepsilon})).
\end{equation}
Now observe that $\det(M-\lambda I)$ is a polynomial of degree $2N$ in $\lambda$, call it $P_\varepsilon(\lambda)$.  Moreover,
\begin{equation}
\lim_{\varepsilon\rightarrow 0}\frac{1}{\varepsilon^N}P_\varepsilon(\sqrt{\varepsilon}\gamma)=\det(\gamma_0^2 I+2 V_{\theta\theta}).
\end{equation}
The latter is a polynomial of degree $N$ in $\zeta=\gamma_0^2$ which has $N-1$ nonzero roots.  Therefore, for $\varepsilon$ sufficiently small there are $2N-2$ zeroes of $P_\varepsilon$ of the form $\lambda=\sqrt{\varepsilon}\gamma(\varepsilon)$ such that $\gamma(\varepsilon)\rightarrow \pm \gamma_0$ and $|\gamma(\varepsilon)|\geq c_1$.  As we remarked above, $M$ has two zero eigenvalues, so those two, together with the $2N-2$ just found comprise the entire spectrum of $M$.\newline
\indent The preceding arguments hold for $\varepsilon\neq0$, regardless of sign.  For the rest of the proof we assume $\varepsilon>0$ and remark that completely analogous arguments can be made for $\varepsilon<0$.  \newline
\indent Suppose that $\varepsilon>0$ and $\phi$ is not a local minimum of the potential.  Then $2V_{\theta\theta}(\phi)$ must have a negative eigenvalue which implies that for all $\varepsilon$ sufficiently small, $\gamma(\varepsilon)$ must have nonzero real part. It follows that $M$ has eigenvalues with nonzero real part for $\varepsilon$ small and so the relative equilibrium is not linearly stable. \newline
\indent The final piece of the proof is to show that in the case that $\varepsilon >0$ and $\phi$ is a local minimum of the potential, the hypothesis of Lemma \ref{lma:3} holds.  Let $v=(v_r,v_\theta)^T$ be an eigenvector of $M$ associated to the eigenvalue $\lambda=\sqrt{\varepsilon}\gamma(\varepsilon)$.  Using the equation $(M-\lambda I)v=0$ we see that
\begin{equation}
-\lambda v_r=\mathcal{O}(\varepsilon)
\end{equation}
which implies that $v_r=\mathcal{O}(\sqrt{\varepsilon})$.  Set $v_r=\sqrt{\varepsilon}w_r(\varepsilon)$ so that the second component of the eigenvector equation becomes
\begin{equation}
-2\sqrt{\varepsilon}w_r-\sqrt{\varepsilon}\gamma v_\theta=\mathcal{O}(\varepsilon).
\end{equation}
From this we see that $w_r=-\frac{\gamma}{2} v_\theta +\mathcal{O}(\sqrt{\varepsilon})$.  Observe also that if we impose the normalization condition $|v_r|^2+|v_\theta|^2=1$ we can conclude that $|v_\theta|^2=1+\mathcal{O}(\varepsilon)$.\newline
\indent Suppose that $\zeta\in\mathbb{R}$ is a nonzero eigenvalue of $2 V_{\theta\theta}$.  It follows from the assumption that $\phi$ is a local minimum that $\zeta>0$, and so $-\gamma^2=\zeta+o(1)$  then implies that $\gamma=\beta i+o(1)$ where $\beta\in\mathbb{R}$ is nonzero.  \newline
\indent  Finally, a short computation gives
\begin{align}
\Omega(v,\bar{v})&=2\textrm{Im}(\bar{v}_r\cdot v_\theta)\notag\\
&=2 \textrm{Im}\left(\frac{\beta i\sqrt{\varepsilon}}{2}|v_\theta|^2 +\mathcal{O}(\varepsilon)\right)\\
&=\beta\sqrt{\varepsilon}+\mathcal{O}(\varepsilon)\notag
\end{align}
where the last equality follows from the normalization assumption.  Therefore, for $\varepsilon$ small and positive the hypothesis of Lemma \ref{lma:3} is satisfied and the desired result follows.
\end{proof}

\begin{remark}\label{rmk:4} According to the theorem, relative equilibria which are near nondegenerate local maxima or minima of the potential undergo a complete loss or gain of stability as the circulation strength $\varepsilon$ passes through zero, while solutions near saddle points remain unstable.  Thus, the behavior of the system if very different depending on whether or not the circulation of the ``small'' vortices has the same sign, or opposite sign, as that of the large central vortex.\end{remark}

\section{Results for Small N}\label{sect:4}

In this section we describe solutions of the $(1+N)$-vortex problem for $N=2,3$ and $4$. For $N=2$ and $3$, we discover all possible relative equilibria of the problem and examine linear stability.  For $N=4$ we use a numerical root finder based on Newton's method to find critical points of the potential and also examine stability.

\subsection{$N=2$}\label{sect:4.1}
It is well known that the only relative equilibria of the 3-vortex problem are collinear and equilateral triangle configurations. First consider the triangle. We compute the Hessian matrix of $V(\theta_1,\theta_2)$ with, e.g., $\theta_1=0,\theta_2=\frac{\pi}{3}$.  Then
\[V_{\theta\theta}(0,\frac{\pi}{3})=
\left( {\begin{array}{cc}
 \frac{3}{2}&-\frac{3}{2} \\
 -\frac{3}{2} & \frac{3}{2}  \\
 \end{array} } \right)
\]
which has eigenvalues $3$ and $0$.  Therefore, for $\varepsilon>0$ Theorems \ref{thm:1} and \ref{thm:2} imply the existence of a linearly stable family of relative equilibria converging to this configuration as $\varepsilon\rightarrow 0$.  Similarly, we discover that the collinear configuration corresponds to a linearly unstable family.

\begin{figure}[H]
 \centering
  \subfloat[]{\label{fig:1a}\includegraphics[scale=0.5]{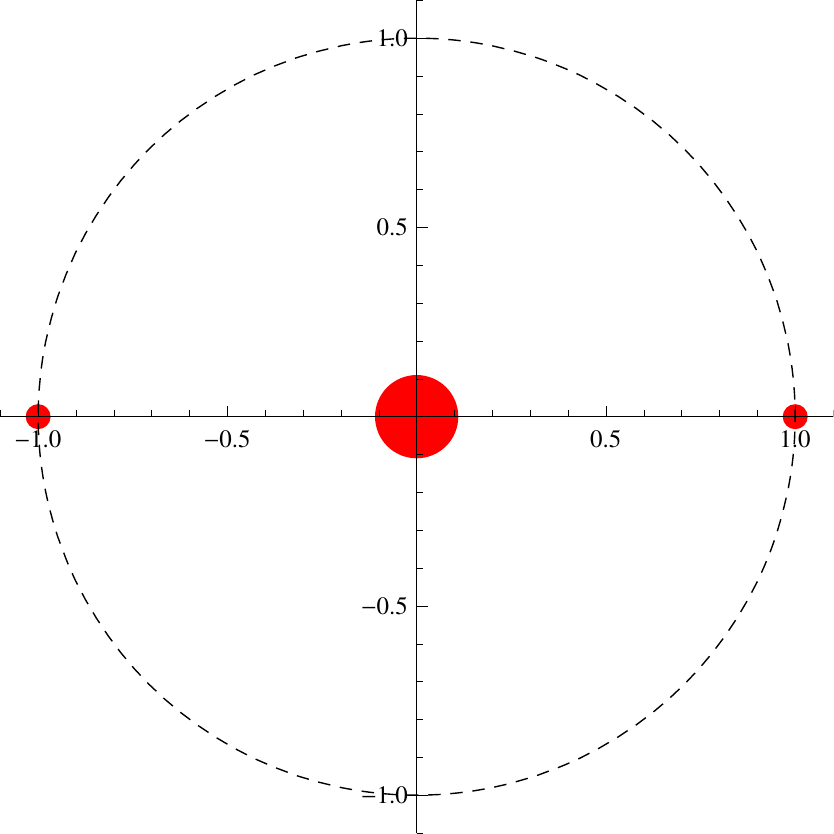}}\;
  \subfloat[]{\label{fig:1b}\includegraphics[scale=0.5]{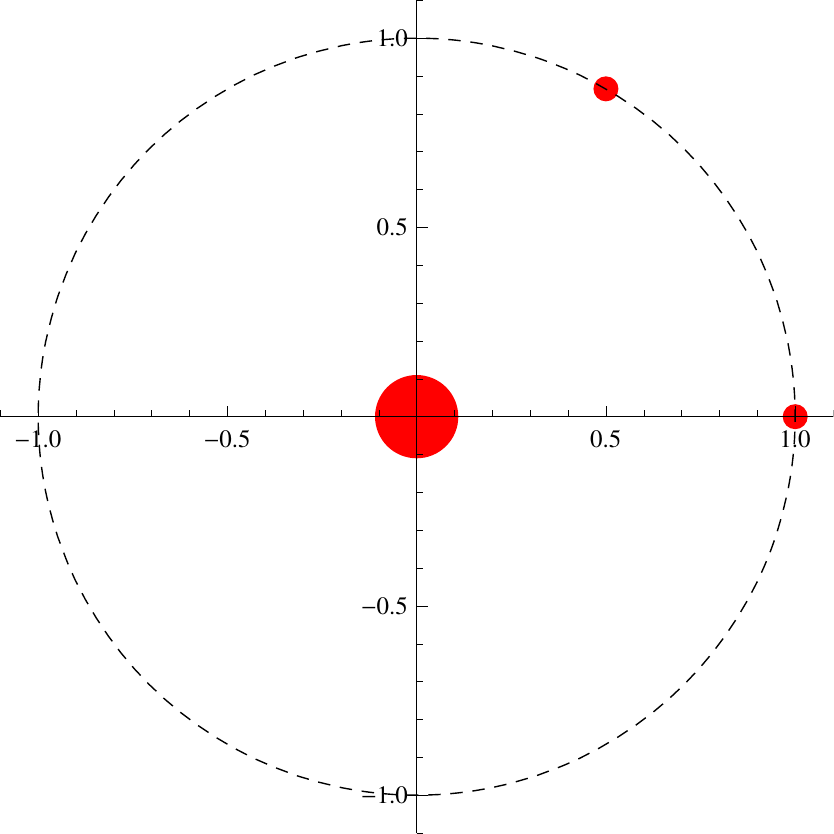}}
    \caption{Relative equilibria of the 1+2 vortex problem.  For $\varepsilon>0$, configuration (a) is the limit of a family of linearly unstable relative equilibria and (b) is the limiting configuration of a sequence of linearly stable equilibria.  For $\varepsilon<0$, the two stability types switch.}
    \label{fig:1}
\end{figure}

\bigskip

\subsection{$N=3$}\label{sect:4.2}
\begin{figure}[H]
 \centering
  \subfloat[]{\label{fig:2a}\includegraphics[scale=0.5]{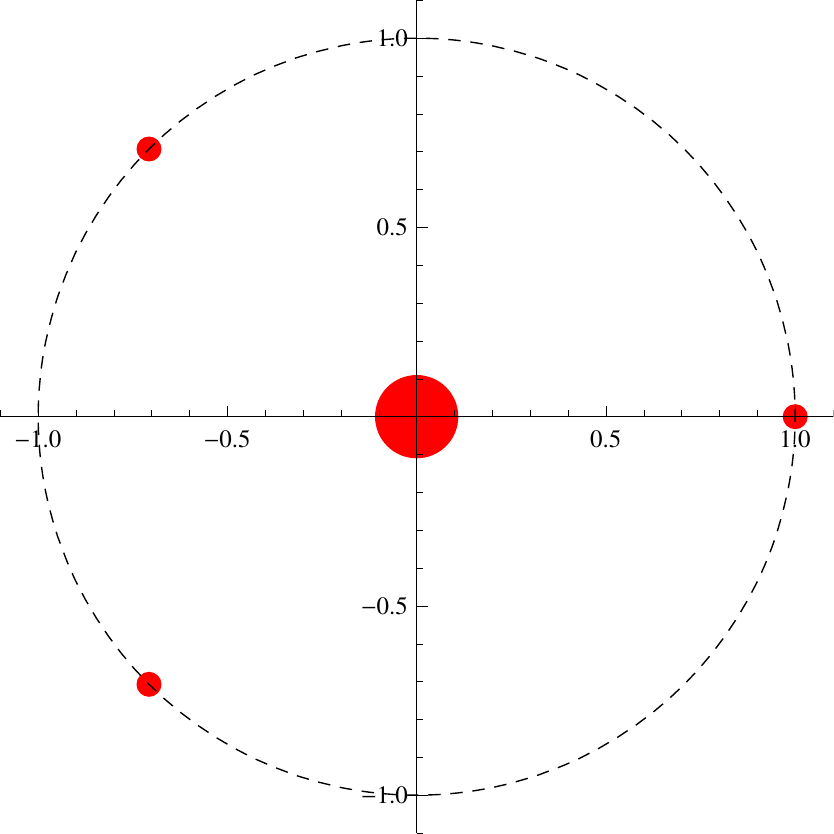}}\;
  \subfloat[]{\label{fig:2b}\includegraphics[scale=0.5]{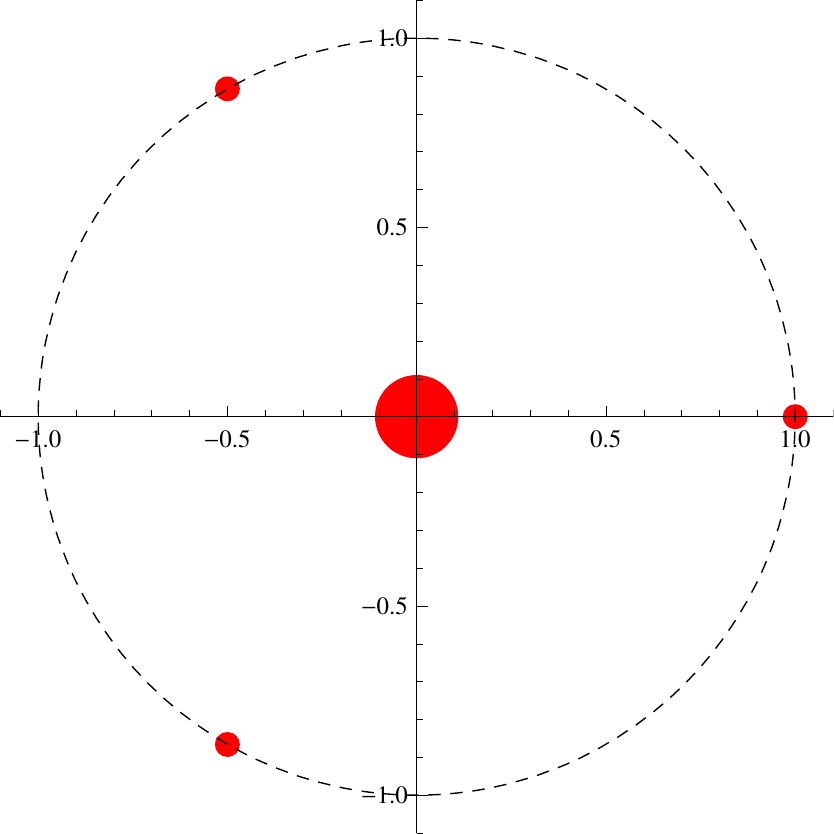}}\;
  \subfloat[]{\label{fig:2c}\includegraphics[scale=0.5]{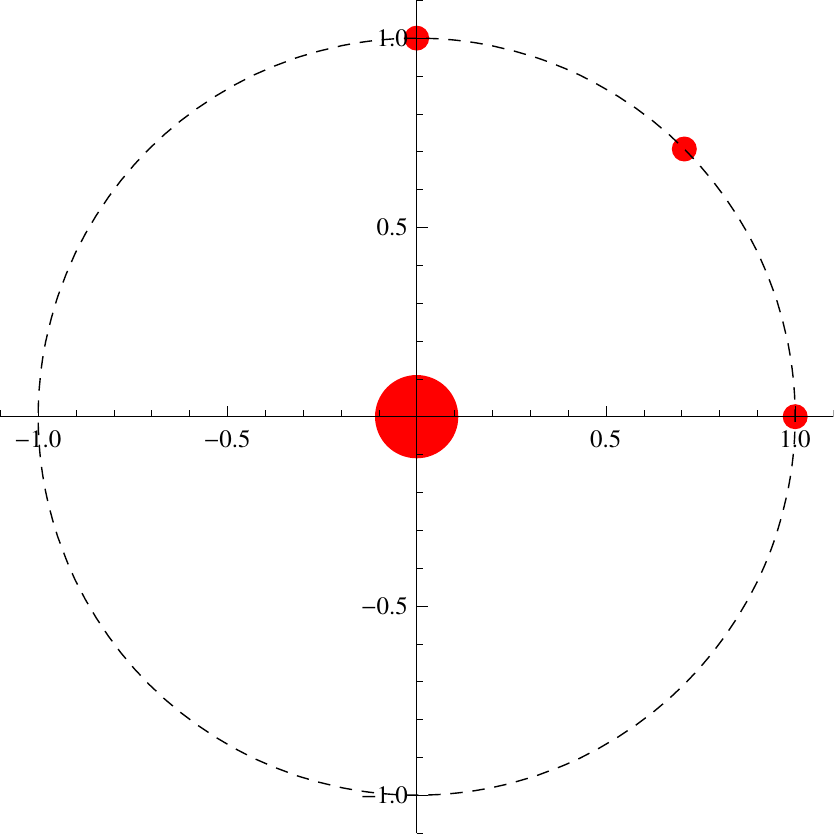}}
  \caption{Relative equilibria of the 1+3 vortex problem. For $\varepsilon>0$ the configurations (a) and (b) are limits of a family of linearly unstable orbits, while configuration (c) corresponds to a family of linearly stable orbits.  For $\varepsilon<0$, configuration (b) is the limit of a family of linearly stable relative equilibria, while the other two correspond to unstable solutions.}
   \label{fig:2}
\end{figure}

By straightforward computations using the potential one can show that the $(1+3)$-vortex problem has exactly three relative equilibria up to rigid rotations and permutations of indices.  Figure \ref{fig:2} illustrates these three families.  In Figure \ref{fig:2a}, the vortices on the left side of the circle each form an angle of $\frac{3\pi}{4}$ with the positive x-axis.  This configuration is a saddle point of the potential and so it corresponds to a linearly unstable family of solutions for both $\varepsilon>0$ and $\varepsilon<0$.  The equilateral configuration about the central vortex is a local maximum of the potential. In Figure \ref{fig:2c}, the small vortices are separated by an angle of $\pi/4$ and this configuration is a local minimum of the potential.

\subsection{$N=4$}\label{sect:4.3}
\begin{figure}[H]
 \centering
  \subfloat[]{\label{fig:3a}\includegraphics[scale=0.5]{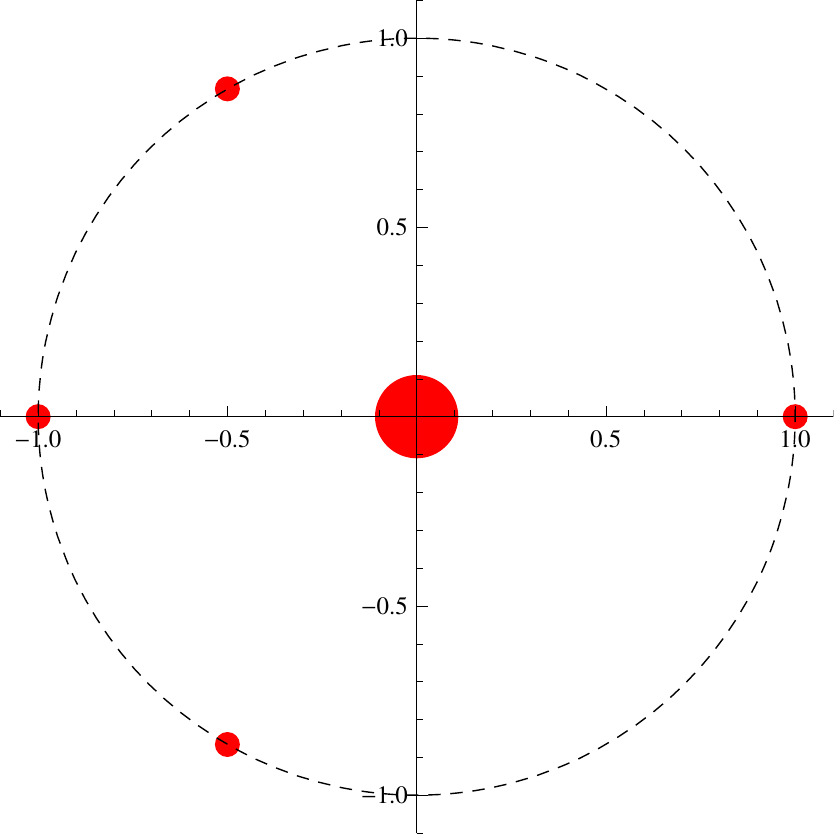}}\;
  \subfloat[]{\label{fig:3b}\includegraphics[scale=0.5]{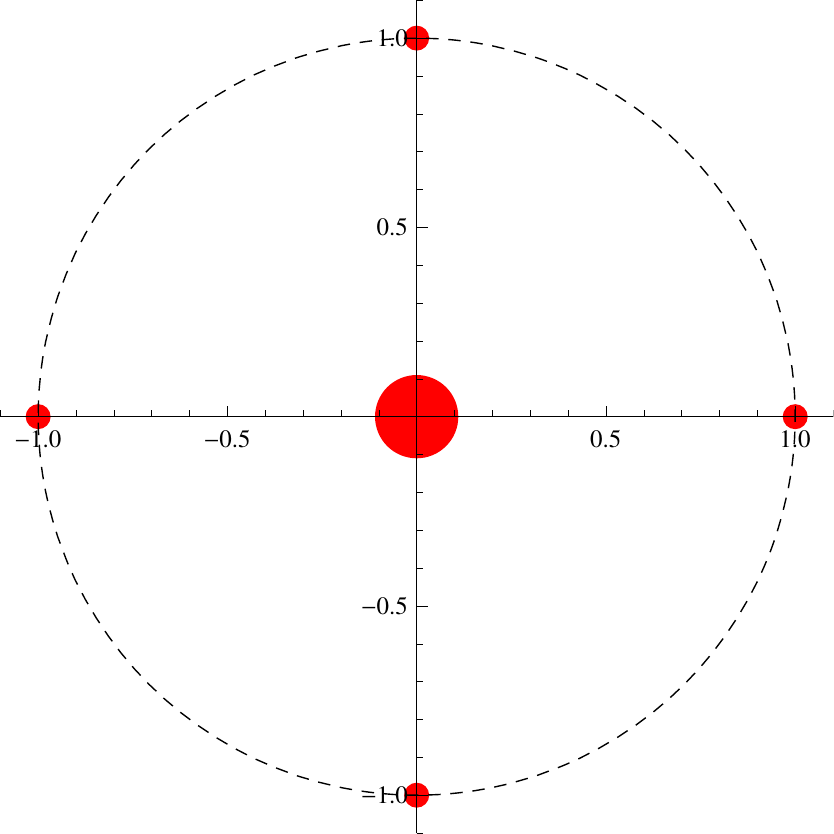}}\;
  \subfloat[]{\label{fig:3c}\includegraphics[scale=0.5]{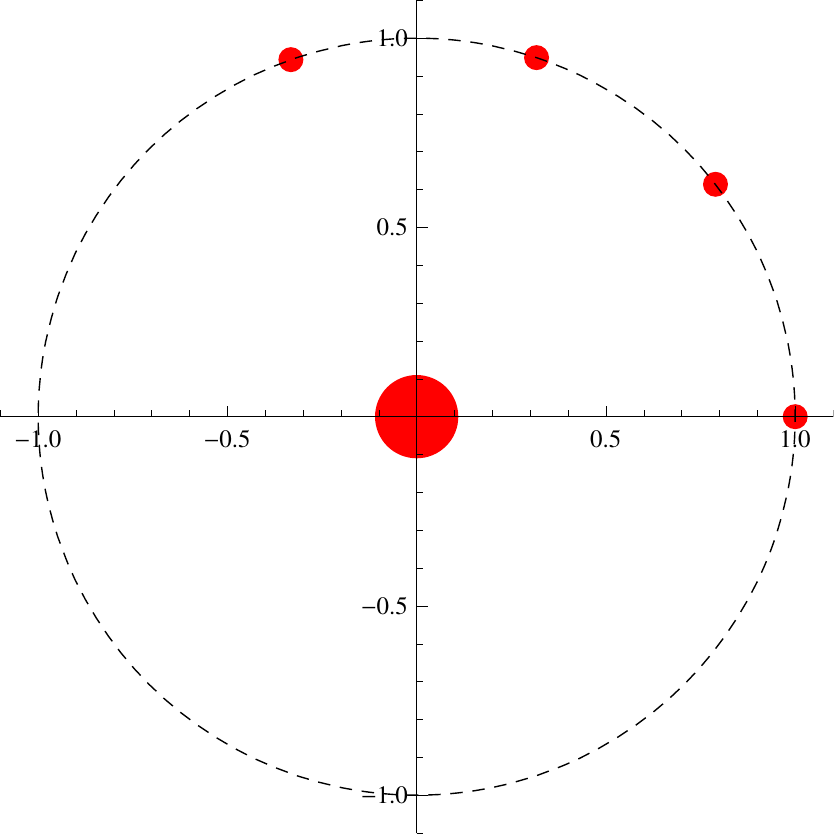}}
  \caption{Numerically observed relative equilibria of the (1+4)-vortex problem. For $\varepsilon>0$, they correspond to an (a) unstable, (b) unstable, (c) stable family of relative equilibria.  For $\varepsilon<0$, all three of the converging families are linearly unstable.}
   \label{fig:3}
\end{figure}

The graphics in Figure \ref{fig:3} show the numerically observed relative equilibria for $N=4$.  In configuration \ref{fig:3a}, two of the vortices form an angle of $\frac{\pi}{3}$ with the negative x-axis.  The remaining two small vortices are at $1$ and $\pi$.  Using the potential, we calculate the nonzero eigenvalues of the Hessian at this configuration to be $4,-\frac{3}{2},1$ and so this is a saddle point.  We next calculate the nonzero eigenvalues associated to the $N$-gon to be $2,-\frac{1}{2},-\frac{1}{2}$.  Thus the $N$-gon is also a saddle point of the potential. In fact, we will show in the next section that the $N$-gon is always a saddle point of the potential for $N\geq 4$.  The vortices in the third configuration are separated by an angle slightly less than $\frac{\pi}{5}$.  The eigenvalues are approximately $12.4, 8.4,$ and $3.7$ and so we find this to be a local minimum of the potential.

We would like to emphasize the fact that all three families of relative equilibria become unstable as $\varepsilon$ passes through zero and becomes negative.  This result also applies for larger values of $N$, which we describe in the next section.

\section{Results for Large N}\label{sect:5}

In this section we prove that the $N$-gon is always a nondegenerate saddle point of the limit potential if $N\geq 4$.  This implies that the limit potential has additional critical points, one of which must be a minimum.  In the simplest case, when all critical points are nondegenerate, there must be at least three distinct families.  In numerical searches, we find just these three families for $N\geq 4$, and they are nondegenerate in every case computed.  When nondegenerate, Theorem 2 implies these critical points continue to $\epsilon\neq 0$. The minima continue to linearly stable equilibria for $\epsilon>0$.  When $\varepsilon<0$, the relative equilibria corresponding to both the $N$-gon and the minimum become highly unstable.  We investigate the hypothesis that this is also the case for the third observed family, thus highlighting the dramatic loss of stability that occurs when the circulation of the ``small'' vortices is opposite in sign from that of the ``large'' central vortex.

\subsection{Instability of the N-gon and Implications}\label{sect:5.1}

For all $N$, the $N$-gon about a central vortex at the origin is a relative equilibrium of the $(1+N)$-vortex problem.  In the previous section we explained that when $N=2$ or $3$ any sequence of relative equilibria converging to the $N$-gon is made up of linearly unstable orbits for $\varepsilon>0$ and linearly stable orbits if $\varepsilon<0$.  When $N=4$ the $N$-gon is a saddle point of the potential and so for both $\varepsilon>0$ and $\varepsilon<0$ it corresponds to families of linearly unstable periodic orbits.  We now show that this last statement generalizes to $N>4$.

\begin{theorem}\label{thm:3} For $N\geq 4$, the Hessian matrix of the potential $V(\theta)$ evaluated at the $N-$gon has at least one positive eigenvalue and one negative eigenvalue.\end{theorem}

\begin{proof}  The Hessian has the form of a circulant matrix with entries
\begin{align}
V_{\theta_i\theta_j}&=-\cos(\theta_i-\theta_j)-\frac{1}{2-2\cos(\theta_i-\theta_j)},\;\; i\neq j\\
V_{\theta_i\theta_i}&=-\sum_{j\neq i}V_{\theta_i\theta_j}.
\end{align}
Let $(a_0,a_1,...,a_{N-1})$ denote the first row of this matrix. The eigenvalues of such a matrix can be computed by evaluating the generator polynomial
\begin{equation}
q(t)=a_0 +a_1 t+c_2 t^2+...+ a_{N-1}t^{N-1}
\end{equation}
of the matrix at $t=\omega^j$, $j=0,...,N-1$ and $\omega=\exp{\frac{2\pi i}{N}}$ is the primitive $N$th root of unity.
Observe first that $q(1)=0$.  Next we evaluate $q(\omega)$:
\begin{align}
q(\omega)&=-\sum_{j\neq 1} V_{\theta_1\theta_j}+\omega V_{\theta_1\theta_2}+\omega^2 V_{\theta_1\theta_3}+...+\omega^{N-1}V_{\theta_1\theta_{N}}\\
&=V_{\theta_1\theta_2}(\omega-1)+V_{\theta_1\theta_3}(\omega^2-1)+...+V_{\theta_1\theta_{N}}(\omega^{N-1}-1)\\
&=\sum_{j\neq 1} V_{\theta_1\theta_j}(\omega^{j-1}-1)\\
&=\sum_{j\neq 1} \left(-\cos(\theta_1-\theta_j)-\frac{1}{2-2\cos(\theta_1-\theta_j)}\right)\left(\exp{\frac{2(j-1)\pi i}{N}}-1\right).
\end{align}
In the case of the $N$-gon, if we assume that the angles are ordered counterclockwise around the unit circle with $\theta_1=0$ then $\theta_j-\theta_1=\frac{2 (j-1) \pi}{N}$.  Further, since this matrix is Hermitian all of its eigenvalues must be real and so we can ignore the imaginary part of the sum.  Therefore, we obtain
\begin{align}
\sum_{k=1}^{N-1}\left(\cos\left(\frac{2 k\pi}{N}\right)+\frac{1}{2-2\cos\left(\frac{2 k\pi}{N}\right)}\right)\left(1-\cos\left(\frac{2 k\pi}{N}\right)\right)=-\frac{1}{2}.
\end{align}
Finally, we compute $q(\omega^2)$ which gives
\begin{equation}
\sum_{k=1}^{N-1} \left(\cos\left(\frac{2k\pi}{N}\right)+\frac{1}{2-2\cos\left(\frac{2k\pi}{N}\right)}\right)\left(1-e^{-\frac{4k\pi i}{N}}\right)=N-2
\end{equation}
whenever $N\geq 4$.\end{proof}

The theorem was proved by computing the first three eigenvalues.  However, it will be useful for what follows to produce estimates on more eigenvalues.  We do this now.  In general, for $0<j<N-1$,
\begin{align}
q(\omega^j)&=\sum_{k=1}^N\left(\cos\left(\frac{2k\pi}{N}\right)+\frac{1}{2-2\cos\left(\frac{2k\pi}{N}\right)}\right)\left(1-\omega^{-jk}\right)\\
&=\sum_{k=1}^N \cos\left(\frac{2k\pi}{N}\right)-\sum_{k=1}^N \cos\left(\frac{2k\pi}{N}\right)\omega^{-jk}+\frac{1}{2}\sum_{k=1}^N \frac{1-\omega^{jk}}{2-2\cos\left(\frac{2k\pi}{N}\right)}.
\end{align}
We examine each of the three sums separately.  Observe first that
\begin{equation}
\sum_{k=1}^N \cos\left(\frac{2k\pi}{N}\right)=-1.
\end{equation}
As in the proof of Theorem \ref{thm:3} we ignore the imaginary parts of the sums.  Next observe that the third term satisfies
\begin{equation}
\frac{1}{2}\sum_{k=1}^N \frac{1-\cos\left(\frac{2jk\pi}{N}\right)}{1-\cos\left(\frac{2k\pi}{N}\right)}=b(j,N)>0
\end{equation}
for $j=2,...,N-2$.  Now write the second term as follows:
\begin{align}
\sum_{k=1}^N \cos\left(\frac{2k\pi}{N}\right)\omega^{-jk}&=\frac{1}{2}\sum_{k=1}^N \left(\omega^{k(1-j)}+\omega^{-k(1+j)}\right).
\end{align}
We can use geometric series to compute this sum, except when $j=1,N-1$ and when $j=\frac{N}{2}+1, \frac{N}{2}-1$ if $N$ is even.  For $j=1$ and $N-1$, the entire sum is equal to $-\frac{1}{2}$.  For all other $j$, we compute the second term, above, to be $-1$, even in the special cases listed.  Therefore, we find that $V_{\theta\theta}$ at the $N$-gon has one zero eigenvalue, two eigenvalues which are $-\frac{1}{2}$ and all other eigenvalues, $\lambda=q(\omega^j)$, satisfy
\begin{equation}
q(\omega^j)=-1-(-1)+b(j,N)=b(j,N)>0.
\end{equation}
Note that since the $N$-gon has only two negative eigenvalues, its continuation has only two directions of instability when $\varepsilon>0$.  However, when $\varepsilon<0$ these correspond to two directions of stability, while the $N-3$ positive eigenvalues correspond to directions of instability.  Thus, for large $N$ solutions near the $N$-gon are much more unstable when $\varepsilon<0$.\newline
\indent We now compare Theorem \ref{thm:3} to the results in \cite{MR1740937}.  In the article, the authors showed that the $N$-gon about a central vortex is linearly, and in fact Liapunov, stable for a bounded interval of possible strengths of the central vortex and that it is unstable whenever the strength is outside of this interval.  In particular, they showed that if $\Gamma$ is the common strength of the vortices lying on the $N$-gon, and $p\Gamma$ is the strength of the central vortex, then the configuration is Liapunov stable if and only if $N$ satisfies

\begin{align}
&\frac{N^2-8N+8}{16}<p<\frac{(N-1)^2}{4}\;\; \text{for N even},\\
&\frac{N^2-8N+7}{16}<p<\frac{(N-1)^2}{4}\;\; \text{for N odd}.
\end{align}

In the setting of this paper, $\Gamma=\varepsilon$ and so $p=\frac{1}{\varepsilon}$.  Thus for $\varepsilon$ sufficiently small (as guaranteed by the Implicit Function Theorem in our arguments), $p$ does not fall into the interval for stability and so the configuration is unstable.  This agrees with our result.  Moreover, it is interesting to note that their result forces an upper bound of $\mathcal{O}\left(\frac{1}{N^2}\right)$ on $\varepsilon$.  See \ref{MR1816912} for related results in the $(1+N)$-body problem.\newline
\indent We finish this section with an important corollary to Theorem \ref{thm:3}.

\begin{corollary}\label{cor:1} Assuming all critical points of $V$ are nondegenerate, there are at least three distinct families of relative equilibria of the $(1+N)$-vortex problem for $N\geq 3$, one of which continues to linearly stable equilibria when $\varepsilon>0$.\end{corollary}

\begin{proof} The case $N=3$ was considered in Section \ref{sect:4.2}.  For $N\geq 4$, we showed that the $N$-gon about a central vortex makes up one of these families, and that it is linearly unstable for $\varepsilon\neq 0$ sufficiently small.  We now show that $V$ has a local minimum.  This family continues to linearly stable relative equilibria when $\varepsilon>0$ by Theorem \ref{thm:2}.  Without loss of generality, assume $\theta_1=0$ and $0<\theta_2<\theta_3<...<\theta_N<2\pi$.  All other critical points of the potential are rotations or permutations of critical points of this form.  Then $V$ is a continuous function on the interior of the specified ``wedge'' of the $N$-cube $(0,2\pi)^N$ that tends to positive infinity at the boundaries (e.g. as $\theta_2\rightarrow \theta_3$).  Thus $V$ must achieve a local minimum inside the wedge.\newline
\indent With restrictions on $\theta$ as above, consider the transformation
\begin{equation}
(\theta_1,\theta_2,...,\theta_N)\mapsto (\theta_1, \theta_2-\theta_1,...,\theta_N-\theta_{N-1})=:(\eta_1,\eta_2,...,\eta_N)
\end{equation}
where $\eta_1= 0$, $\eta_i>0$, $i=2,...,N$ and $\eta_1+\eta_2+...+\eta_N<2\pi$.  Since $V(\eta)$ does not depend explicitly on $\eta_1$ we simply drop this coordinate and consider the function of $N-2$ variables $V(\eta_2,...,\eta_{N-1})$.  This removes the zero eigenvalue from the associated Hessian.
Under this transformation, the signs of the remaining eigenvalues of the Hessian are preserved.  Then in the case of the $N$-gon, $V_{\eta\eta}$ has exactly 2 negative eigenvalues.  Therefore, the index of the $N$-gon for the gradient vector field $V_\eta$ is one.  The index of the minimum is also one.  Thus we may apply the Hopf Index Theorem (see, for instance, the standard text \cite{MR658304}) to $V_\eta$ and obtain a third critical point of $V(\eta)$ having index negative one.  This implies the existence of a third critical point in the original coordinates which is a relative equilibrium of the $(1+N)$-vortex problem and is distinct from the $N$-gon and the minimum. \end{proof}

We remark that if the minimum is isolated then it can be continued by topological arguments to the $\varepsilon\neq 0$ case.  While more degenerate situations are conceivable, all the critical points seen numerically have been nondegenerate.  Assuming this is true of the third family given by Corollary \ref{cor:1}, note that the associated Hessian matrix must have an odd number of negative eigenvalues, and so it must be a saddle point of the potential when $N$ is odd.  In this case, Theorem \ref{thm:2} implies that its continuation to a relative equilibrium of the full problem is linearly unstable for all $\varepsilon\neq 0$ sufficiently small.  Numerics suggest that for all $N$, even or odd, this third critical point is nondegenerate with exactly one negative eigenvalue.  For this reason we expect its continuation to be unstable for $\varepsilon>0$ and highly unstable for $\varepsilon<0$.  Thus, there is evidence that all three families of relative equilibria become very unstable as $\varepsilon$ becomes negative.\newline

\subsection{Numerical Observations}

\begin{figure}[H]
 \centering
  \subfloat[]{\includegraphics[scale=0.5]{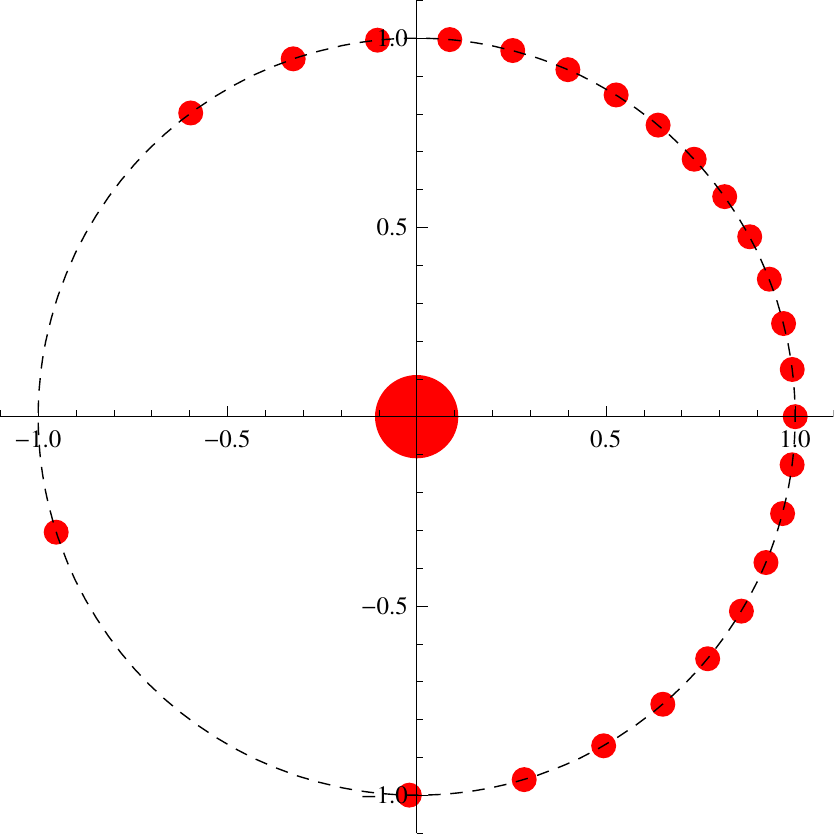}}\;
    \subfloat[]{\includegraphics[scale=0.5]{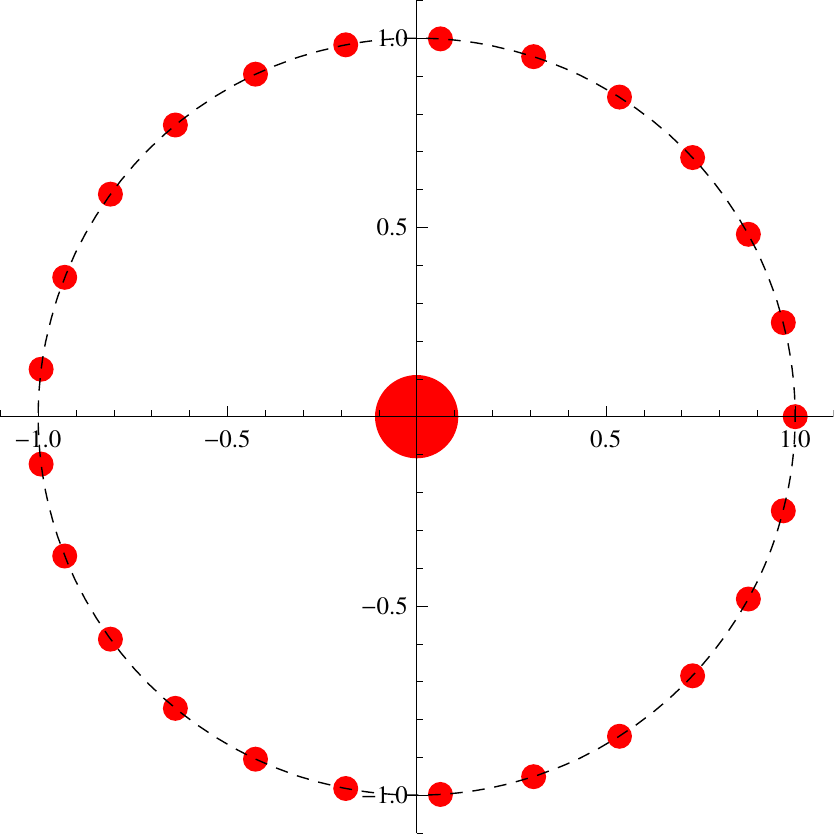}}\;
        \subfloat[]{\includegraphics[scale=0.5]{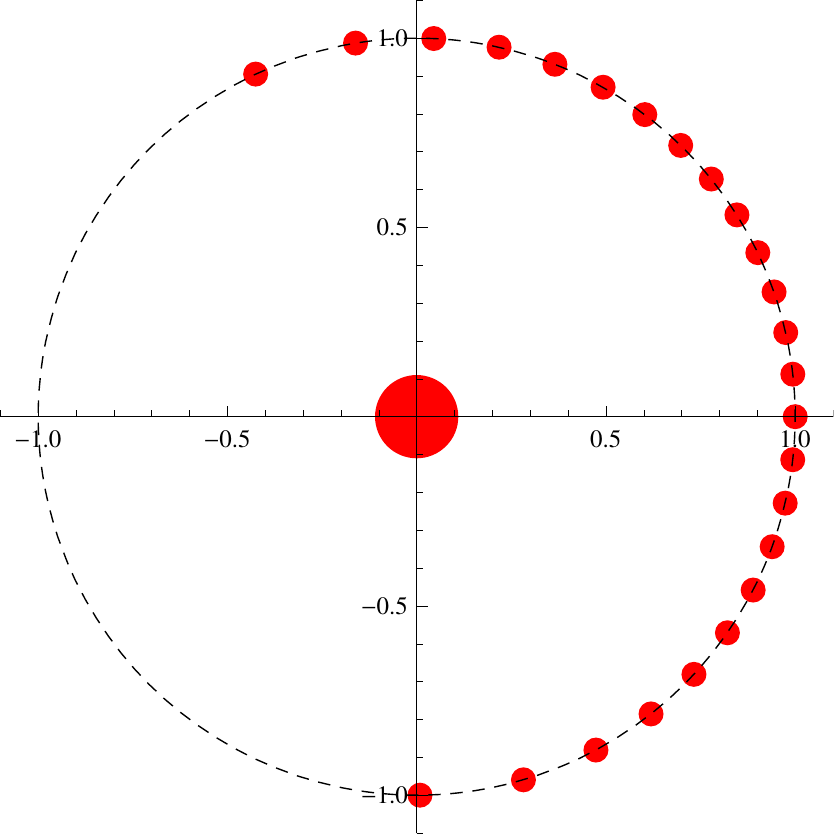}}
    \caption{Three relative equilibria of the $(1+25)$-vortex problem.}
    \label{fig:4}
\end{figure}

We supplement the discussion in the previous subsection with numerical findings.  Using an algorithm based on Newton's method, we were able to locate three distinct families of critical points of the potential.  We have run this algorithm for $N\leq 100$, and in all of these cases the critical points are relative equilibria of the $(1+N)$-vortex problem.  In other words, they continue to relative equilibria of the full problem with $\varepsilon\neq 0$ sufficiently small.\newline
\indent A representative from each of the three families is shown in Figure \ref{fig:4} for $N=25$.  Compare these equilibrium types to those in Figures \ref{fig:1}, \ref{fig:2} and \ref{fig:3}.  It is natural to wonder if these three critical points correspond to relative equilibria of the full problem for all $N$, and if they are exactly those critical points guaranteed by Corollary \ref{cor:1}.  We cannot answer this question conclusively at this time, but we have numerically observed only these three families for $N\leq 100$.\newline
\indent For the two relative equilibria with ``clusters'' of small vortices, we have observed that as $N$ increases the cluster tends to fill out the unit circle more and more, and appears to approach the $N$-gon.  This filling in process is beginning to become apparent in Figure \ref{fig:4}.

\begin{figure}[H]
 \centering
  \subfloat[]{\includegraphics[scale=0.5]{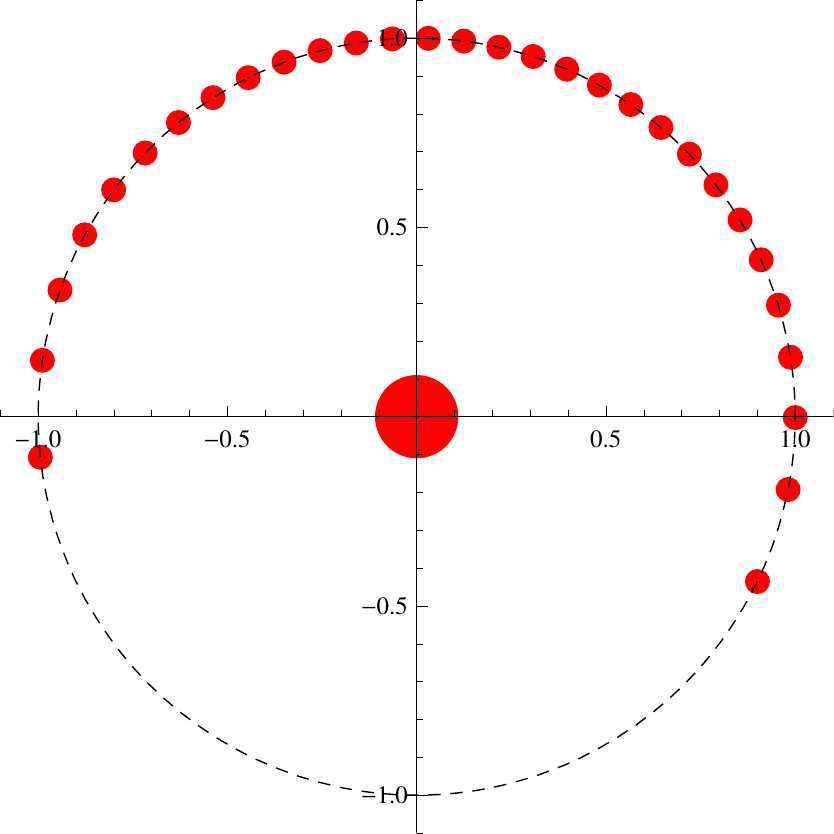}}\;
    \subfloat[]{\includegraphics[scale=0.5]{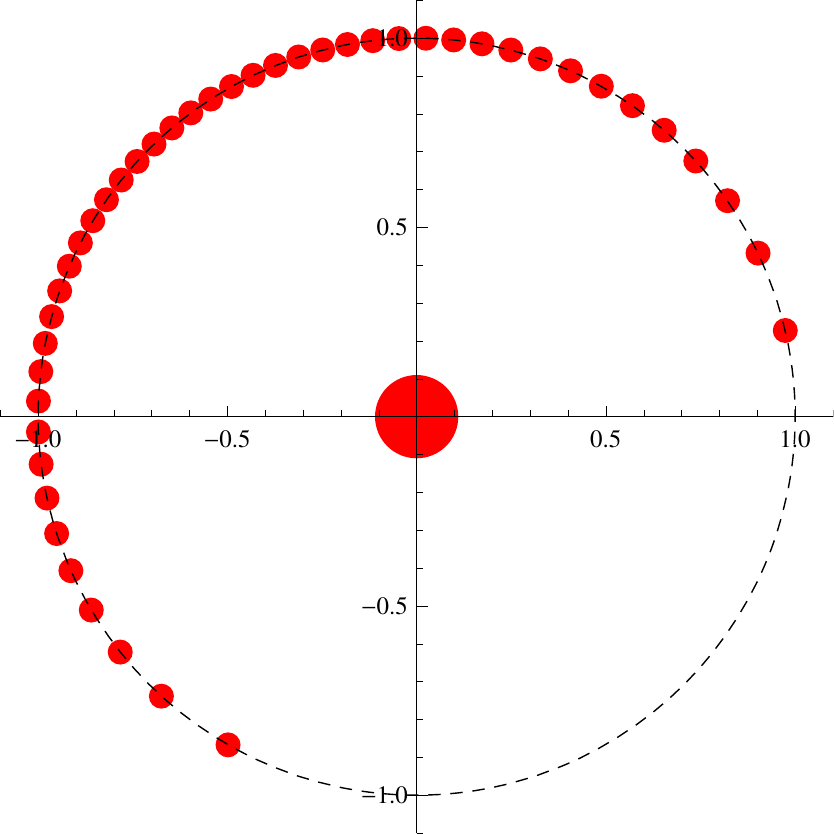}}\;
        \subfloat[]{\includegraphics[scale=0.5]{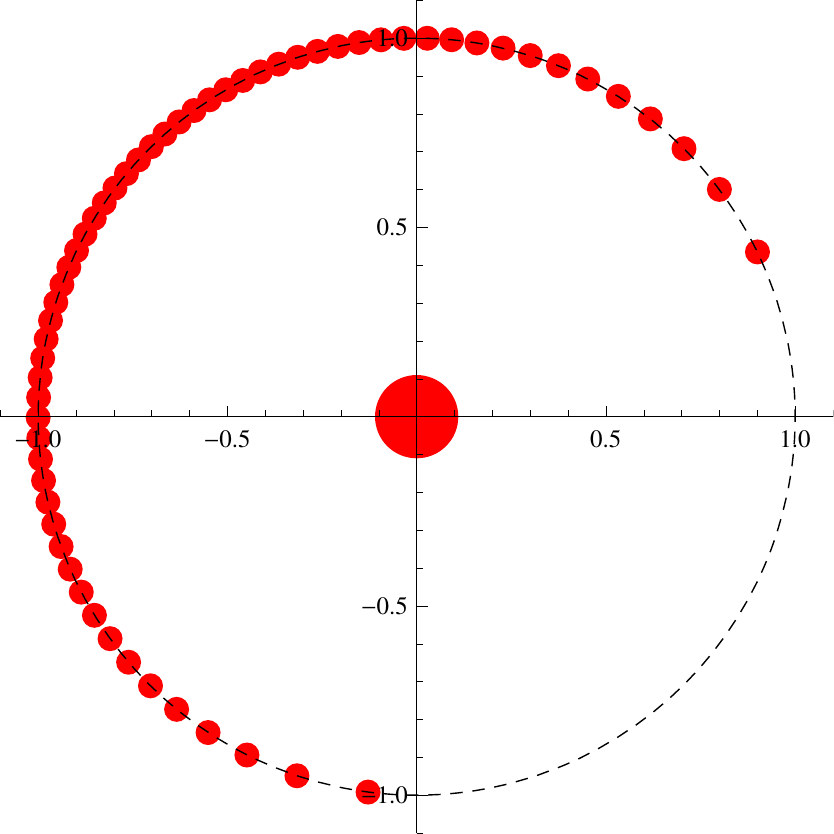}}
    \caption{Some relative equilibria of the $(1+N)$-vortex problem which are minima of the potential. (a) $N=30$, (b) $N=45$, (c) $N=60$.}
    \label{fig:5}
\end{figure}

\indent Using a modification of the algorithm, we were able to identify one of the three classes as containing nondegenerate minima of the potential.  Therefore this class corresponds to linearly stable sequences of relative equilibria of the full problem when $\varepsilon>0$.  Three representatives are shown in Figure \ref{fig:5}.

\indent To investigate the stability of the third class of critical points we computed eigenvalues of the associated Hessian matrix for $N\leq 100$.  We discovered that the matrix has exactly one negative eigenvalue and $N-2$ positive eigenvalues, thus leading us to believe that, like the $N$-gon, nearby relative equilibria are highly unstable when $\varepsilon<0$.  We reiterate that this third family of relative equilibria has an odd number of negative eigenvalues (at least for $N\leq 100$) and is therefore a candidate for the third family of critical points guaranteed by the Hopf Index Theorem in the proof of Corollary \ref{cor:1}.\newline
\indent Finally, we remark that for all three observed families of critical points, the positive eigenvalues of the Hessian appear to increase monotonically with $N$.  This is still more evidence for the striking instability of all three relative equilibria of the full problem when $\varepsilon<0$.

\section{Conclusions}

In this article, we examined solutions of a specialized problem in point vortex dynamics: the $(1+N)$-vortex problem.  We were able to show that all solutions of this problem must be critical points of a specified potential function and that whenever such a critical point is nondegenerate in some appropriate sense, there exists a nondegenerate family of relative equilibria of the full $(N+1)$-vortex problem which converges to it.  Moreover, we were able to exploit properties of Hamiltonian systems to show that the linear stability of a member of such a family is intimately tied to the potential.\newline
\indent We used these results to prove that the $N$-gon about a central vortex is a linearly unstable relative equilibrium when the ``weak'' vortex circulation is sufficiently small if $N\geq 4$.  We further showed that when all critical points of the limit potential are nondegenerate, there are at least two other distinct families of relative equilibria with $N$ vortices having small, nonzero circulation.  One of these families was shown to be linearly stable when the central and surrounding vortices are of the same sign.  Our further numerical investigation supports the hypothesis that the problem has no linearly stable relative equilibria when the weak and strong circulations are opposite in sign.

\bigskip

\noindent \textbf{Acknowledgement:} The research of A.M.B. was supported by the Center for BioDynamics at Boston University and the NSF (DMS$0602204$ EMSW$21$-RTG). C.E.W. was supported in part by DMS-$0908093$.  He would also like to thank Paul Newton for a very useful discussion of relative equilibria of point vortices.

\appendix
\section{Appendix}
The following lemma was used in the proof of Theorem \ref{thm:2}.

\begin{lemma}\label{lma:4} Let $A$, $B$, $C$, and $D$ be real $n\times n$ matrices and define
\begin{equation}
M:=\left( \begin{array}{cc}
A & B  \\
C & D
\end{array} \right).
\end{equation}
Then $\det M=\det(A)\det(D-CA^{-1}B)=\det(D)\det(A-BD^{-1}C)$ whenever $A$ and $D$ are invertible.\end{lemma}

\begin{proof} The result follows immediately from the following observation:
\begin{align}
M&=\left( \begin{array}{cc}
A & 0  \\
C & I
\end{array} \right)
\left( \begin{array}{cc}
I & A^{-1}B  \\
0 & D-CA^{-1}B
\end{array} \right)\\
&=\left( \begin{array}{cc}
I & B  \\
0 & D
\end{array} \right)
\left( \begin{array}{cc}
A-BD^{-1}C & 0  \\
D^{-1}C & I
\end{array} \right).
\end{align}
\end{proof}

\newpage

\def\cprime{$'$}


\begin{thebibliography}{10}

\bibitem{MR1332658}
Hassan Aref.
\newblock On the equilibrium and stability of a row of point vortices.
\newblock {\em J. Fluid Mech.}, 290:167--181, 1995.

\bibitem{MR2337012}
Hassan Aref.
\newblock Point vortex dynamics: a classical mathematics playground.
\newblock {\em J. Math. Phys.}, 48(6):065401, 23, 2007.

\bibitem{MR658304}
Raoul Bott and Loring~W. Tu.
\newblock {\em Differential forms in algebraic topology}, volume~82 of {\em
  Graduate Texts in Mathematics}.
\newblock Springer-Verlag, New York, 1982.

\bibitem{MR1740937}
H.~E. Cabral and D.~S. Schmidt.
\newblock Stability of relative equilibria in the problem of {$N+1$} vortices.
\newblock {\em SIAM J. Math. Anal.}, 31(2):231--250 (electronic), 1999/00.

\bibitem{MR1310189}
Josefina Casasayas, Jaume Llibre, and Ana Nunes.
\newblock Central configurations of the planar {$1+n$} body problem.
\newblock {\em Celestial Mech. Dynam. Astronom.}, 60(2):273--288, 1994.

\bibitem{MR2104897}
Josep~M. Cors, Jaume Llibre, and Merc{\`e} Oll{\'e}.
\newblock Central configurations of the planar coorbital satellite problem.
\newblock {\em Celestial Mech. Dynam. Astronom.}, 89(4):319--342, 2004.

\bibitem{hall}
G.R. Hall.
\newblock Central configurations in the planar 1+n body problem.
\newblock 1988.
\newblock Boston University, preprint.

\bibitem{MR987772}
Jair Koiller and Sonia~P. Carvalho.
\newblock Nonintegrability of the {$4$}-vortex system: analytical proof.
\newblock {\em Comm. Math. Phys.}, 120(4):643--652, 1989.

\bibitem{MR1411341}
D.~Lewis and T.~Ratiu.
\newblock Rotating {$n$}-gon/{$kn$}-gon vortex configurations.
\newblock {\em J. Nonlinear Sci.}, 6(5):385--414, 1996.

\bibitem{MR1262722}
Richard Moeckel.
\newblock Linear stability of relative equilibria with a dominant mass.
\newblock {\em J. Dynam. Differential Equations}, 6(1):37--51, 1994.

\bibitem{MR1472896}
Richard Moeckel.
\newblock Relative equilibria with clusters of small masses.
\newblock {\em J. Dynam. Differential Equations}, 9(4):507--533, 1997.

\bibitem{MR1831715}
Paul~K. Newton.
\newblock {\em The {$N$}-vortex problem}, volume 145 of {\em Applied
  Mathematical Sciences}.
\newblock Springer-Verlag, New York, 2001.
\newblock Analytical techniques.

\bibitem{MR648066}
Julian~I. Palmore.
\newblock Relative equilibria of vortices in two dimensions.
\newblock {\em Proc. Nat. Acad. Sci. U.S.A.}, 79(2):716--718, 1982.

\bibitem{MR1816912}
Gareth~E. Roberts.
\newblock Linear stability in the {$1+n$}-gon relative equilibrium.
\newblock In {\em Hamiltonian systems and celestial mechanics ({P}\'atzcuaro,
  1998)}, volume~6 of {\em World Sci. Monogr. Ser. Math.}, pages 303--330.
  World Sci. Publ., River Edge, NJ, 2000.

\bibitem{MR1753026}
Dieter~S. Schmidt.
\newblock Spectral stability of relative equilibria in the {$N+1$} body
  problem.
\newblock In {\em New trends for {H}amiltonian systems and celestial mechanics
  ({C}ocoyoc, 1994)}, volume~8 of {\em Adv. Ser. Nonlinear Dynam.}, pages
  321--341. World Sci. Publ., River Edge, NJ, 1996.

\bibitem{MR564329}
S.~L. Ziglin.
\newblock Nonintegrability of the problem of the motion of four point vortices.
\newblock {\em Dokl. Akad. Nauk SSSR}, 250(6):1296--1300, 1980.

\end{thebibliography}
\end{document}